\documentstyle[twoside,amssymb,12pt]{article}
\newtheorem{prop}{Proposition}[section]

\newtheorem{theo}[prop]{Theorem}

\title{\sc
New examples of  conservative  systems on $S^2$ 
possessing an integral cubic in momenta
} 

\author{{\sc Elena N. Selivanova}\footnote{Supported by DAAD.} \\ }
\begin{document}
\date{}

\maketitle

\thispagestyle{empty}

\section{\bf Introduction}
\noindent
 It has been proved   that on $2$-dimensional  orientable compact 
manifolds of genus
$g>1$ there is no integrable geodesic flow with
an integral  polynomial in momenta, see \cite{kol,koz}.
There is  a conjecture  that all integrable geodesic flows
on $T^2$ possess an integral  quadratic in momenta (for more details see
\cite{BKF}).  All geodesic flows on $S^2$   and $T^2$ 
possessing integrals linear and  
quadratic  in momenta have been described in \cite{kol, BaNe}, using some ideas  going back to
Darboux \cite{Dar}. So far there has been known  only one example of  conservative system on 
$S^2$   possessing an integral  cubic in momenta: the case of 
Goryachev-Chaplygin in the dynamics of a rigid body, see \cite{BKF}. 
\\
The aim of this paper is to propose a new one-parameter family of examples
 of complete integrable conservative systems on $S^2$ 
 possessing an integral  cubic in momenta.
We show that our family  does not include the case of Goryachev-Chaplygin.
\\
The paper is organized as follows:
\medbreak
\noindent
We start with the  investigation of  the following initial
 value problem
$$u'''u'r^4=-7u''u'r^3-2u''^2r^4+u''ur^2-2u'^2r^2+u'ur+u^2,$$ 
\begin{equation}
u(1)=0, u'(1)=1, u''(1)={\tau}-1,
\label{u}
\end{equation}
where $u: {\bf R}^+ \mapsto {\bf R}$, $r \mapsto u(r)$, 
and prove that there exists a  positive constant $T$ 
such that the one-parameter family   
$\{ \Psi_{\tau } \}_{\tau\in (0,T) }$ of solutions  of this problem 
defines a one-parameter family  of smooth conservative systems on $S^2$ 
with  energy
\begin{equation}
 H_{\tau}(r,\varphi,dr,d\varphi)=\frac{r^2d\varphi^2 + dr^2}
 {r^4\Psi_{\tau }'^2(r)} - \left(\Psi_{\tau }''r^2 + \Psi_{\tau }'r - \Psi_{\tau }\right)
\Psi_{\tau }'^2r^2\cos \varphi 
\label{ms}
\end{equation}
in polar coordinates $r$ and $\varphi$,  possessing nontrivial additional integrals cubic in momenta.  
\\
We note that the limit point of our family for ${\tau}=0$  corresponds 
to the function  $\Psi_0=\frac{1}{2}\left(r - \frac{1}{r} \right)$ which defines by 
(\ref{ms}) the standard metric of  constant curvature $K=1$ on $S^2$.   
\section{\bf Existence of $C^{\infty}$ solutions } 
\noindent
In this chapter the problem of the existence of the 
smooth solutions of (\ref{u}) on $(0, +\infty)$
will be considered. 
\\
Initial value problem (\ref{u}) can be replaced
by the initial value problem  
\begin{equation}
x'x'''=xx'' - 2x''^2 + x'^2 + x^2, \  x(0)=0, x'(0)=1, x''(0)={\tau}.
\label{x}
\end{equation}
This reduction is accomplished by the change of 
variable $r\to t=\log r$ which changes (\ref{u})
to (\ref{x}). Initial value problem (\ref{x})
is autonomous and more convenient for further
investigations. For this reason, the theorem which
follows will be stated for (\ref{u}), but
proved for (\ref{x}).
\medbreak
\begin{theo}
There is a positive constant $T$ 
such that solutions of (\ref{u}) exist on $(0, +\infty)$
for $\tau \in (-T, T)$.
\label{ex-th}
\end{theo} 
\noindent 
{\em Proof.} Initial value problem (\ref{x})
has a unique solution $x(t)=\Theta_{\tau}(t)$
which is positive on $(0,\varepsilon)$ and negative
on $(-\varepsilon, 0)$ for a sufficiently small 
$\varepsilon$.
\medbreak
\noindent
Let us consider the case $t>0$. 
\\
The differential equation from  (\ref{x})
can be replaced by the following differential
equation of the second order 
$$\ddot{q}=\frac{1}{q}\left ( 1+2q^2-3q^4+\dot{q} -
7q^2\dot{q}-2{\dot{q}}^2\right) $$
with function $q(t)=R'(t)$ where $R(t)=\log x(t)$.
We may rewrite this equation as a system of differential equations of the first order:
\begin{equation}
\dot q=p, \ \
\dot p=\frac{1}{q}\left ( 1+2q^2-3q^4+p-7q^2p-2p^2 \right).
\label{syst1}
\end{equation}
Since  system (\ref{syst1}) is symmetric with respect to $q\mapsto -q$, 
$t\mapsto -t$,
 it  suffices to consider the case $q>0$.
\\
In order to obtain the phase portrait of (\ref{syst1}) we may 
consider the following smooth system
\begin{equation}
\dot q=qp, \ \ 
\dot p=1+2q^2-3q^4+p-7q^2p-2p^2. 
\label{sms}
\end{equation}
The solutions of (\ref{syst1}) are obtained from the solutions of
 (\ref{sms}) by a
reparametrization. 
System  (\ref{sms}) has four singular points: two saddle points 
$p=1, q=0$, $p=-\frac{1}{2}, q=0$
 and two knots $p=0, q=\pm 1$. 
\medbreak
The aim of our further investigations is to show
that the orbits of (\ref{syst1}), corresponding to
the solutions of (\ref{x}), when $\tau$ belongs to a certain interval,
converge to the singular point $q=1, p=0$.
\medbreak
For any solution $x(t)=\Theta_{\tau}(t)$, $t\in (0, 
\varepsilon)$ of 
(\ref{x}) there is an orbit $\{ \Gamma_{\tau}: p=p_{\tau}(q)\}$ 
in the phase spase of (\ref{syst1}).
\\
For any solution of (\ref{x}) it holds that 
$R(t)\to -\infty$ as 
$t\to 0+$ and, therefore, $q(t)=R'(t)\to + \infty$,
 $p(t)=q'(t)\to -\infty$ as $t\to 0+$.
Thus $p_{\tau}(q)\to -\infty$ as $q\to +\infty$.
\\
We may now consider only the orbits of (\ref{syst1}) where 
$p\to -\infty$ as $q\to +\infty$.
\\
Show that if $\tau_1>\tau_2$ then $p_{\tau_1}(q)>p_{\tau_2}(q)$.
Indeed, for a solution $x(t)=\Theta_{\tau_1}(t)$ of (\ref{x}) it holds
\begin{equation}
\tau_1=\lim_{t\to 0+}{ \frac{\Theta''_{\tau_1}(t)}{\Theta'_{\tau_1}(t)} }=
\lim_{q\to +\infty}{\frac{q^2 + p_{\tau_1}(q)}{q}}.
\label{ta_u}
\end{equation}
Note that the function $x(t)=\sinh t$ satisfies (\ref{x}) when $\tau=0$.
The related orbit of (\ref{syst1}) has then the following form 
$\{\Gamma_0: p=p_0(q)=-q^2 + 1\}$.
\\
So, the orbits of (\ref{syst1}), corresponding to the solutions of (\ref{x}),
for $\tau \ge 0$ converge to the singular point $q=1, p=0$. 
\medbreak
We prove below that orbit $(*)$, where 
$p\to -\frac{1}{2}$ as $q\to 0+$,
 corresponds to the solution
of (\ref{x}) when $\tau$ is equal to a negative constant $-T$ and all
orbits of (\ref{syst1}), lying between $(*)$ and $\Gamma_0$ correspond
to the solutions of (\ref{x}) for $\tau\in (-T,0)$.
\medbreak
Assume first that there exists a constant $\tau_0<0$ 
such that  orbit $\Gamma_{\tau_0}$ does not converge
to the point $q=1, p=0$.
\\
Consider the set $W$ of the orbits of (\ref{syst1}), 
lying between  $\Gamma_{\tau_0}$ and $\Gamma_0$.
Show that for any orbit in $W$ the value $q$ becomes infinite in a 
finite time interval.
\\
In fact, for any solution of (\ref{syst1})  holds  
\begin{equation}
\int^{q(t)}\frac{dq}{p}=t + const.
\label{oc}
\end{equation}
For orbits in $W$ we have $p<-q^2 +1$. Hence, the left hand side of (\ref{oc})
is bounded for any orbit in $W$ as $q\to +\infty$.
Without loss of generality we assume that for these 
 solutions of (\ref{syst1})  $t$ vanishes  as $q$ becomes infinite. 
\\
We conclude that  for any orbit from $W$ it holds by (\ref{ta_u}) 
$$(\tau_0+o(1))q<p+q^2<1$$
and for a corresponding solution of (\ref{x}) we get for $t \to 0+$
$$\tau_0+o(1)<\frac{x''(t)}{x'(t)}<1.$$
Thus, for a corresponding solution of (\ref{x}) the function  $(\log x'(t))'$ is 
bounded in an interval, containing $0+$.
 Hence,  $\lim_{t\to 0+}{x'(t)}$
is bounded  for a solution $x(t)$ of (\ref{x}), 
corresponding to an orbit in $W$ and, therefore, $x(t)=\frac{x'(t)}{q}\to 0$
as  $t\to 0+$. 
\\
Let us consider orbit $(*)$  where $p=p^*(q)$.
If $(*)$ is the same as $\Gamma_{\tau_0}$ then $T=-\tau_0$. 
\\
If $(*)$ is not the same as $\Gamma_{\tau_0}$ then it belongs to $W$ and,
hence, corresponds to the solution of (\ref{x})
for $\tau=-T$, where $T$ is equal to a positive constant  and, moreover, 
$T<-\tau_0$. 
\\
Thus we have shown that for any $\tau\in (-T, T)$ where $T$ is a
 positive constant sytem (\ref{x}) has a solution for all $t\ge 0$.
\medbreak
Let us assume now that there is no such orbit $\Gamma_{\tau_0}$ and, so, 
all orbits of (\ref{syst1}), corresponding to the solutions of (\ref{x}),
 converge to the singular point $q=1, p=0$.
\\
Clearly, orbit $(*)$ corresponds to a solution $x(t)=\eta(t)$ of 
the differential equation in (\ref{x}).
As  mentioned above, the function $q(t)=\frac{\eta '(t)}{\eta(t)}$ 
becomes infinite in a finite  time interval of $t$, say, as $t\to 0+$.
Since $p^*(q)<-q^2 - kq$ for any $k$ and
 for  sufficiently large $q$, for $\eta(t)$ it holds
\begin{equation}
\lim_{t\to 0+}{\frac{{\eta''(t)}}{{\eta'(t)}}}=-\infty.
\label{p}
\end{equation}
Note that $\eta(0)$ is finite because $R(0)-R(t)=\int_{t}^{0}{q(t)dt}<0$ for $t>0$ and 
$\eta(t)=\exp R(t)$. 
So, there are two cases: $\eta(0)=0$ or $\eta(0)\ne 0$.
\\
Assume  that $\eta(0)\ne 0$. It follows that $\eta'(0+)=+\infty$ 
and from (\ref{p}) we get $\eta''(0+)=-\infty$.
\\
Rewrite now the differential equation in (\ref{x})  in the following form 
$$\frac{x'''}{x}=-3\frac{x''-x}{x'} -2\frac{(x''-x)^2}{xx'} + \frac{x'}{x}.$$
Therefore,
it holds $$\frac{\eta'''}{\eta}=-3\frac{p^*(q) +q^2 -1}{q} -2\frac{(p^*(q) + q^2 -1)^2}{q} + q.$$
We obtain   $\eta'''(0+)=-\infty$. This is a contradiction  $\eta''(0+)=-\infty$.
\\
Assume  that  $\eta(0)=0$ and  $\eta'(0)\ne 0$.
 From (\ref {p}) we get $\eta''(0+)=-\infty$. 
Rewriting  the differential equation in (\ref{x})  in the following form
$$x'''=x'-3\frac{(x''-x)x}{x'} - 2\frac{(x''-x)^2}{x'},$$
we obtain again  $\eta'''(0+)=-\infty$.
\\
Therefore, we only have to consider one case  $\eta(0)=0$ and  $\eta'(0)=0$.
Taking into account that ${\eta'}(t)=\eta(t)q(t)$ and $q(t)>0$ as $t\to 0+$, we 
conclude that ${\eta'}(t)\to 0+$ as $t\to 0+$ but on the other hand  from (\ref {p})
it follows that ${\eta''(t)}<0$.
\\
These contradictions finally show  that there is an orbit $\Gamma_{\tau_0}$
which does not converge to the point $q=1, p=0$
and, therefore, for any $\tau\in (-T, T)$, where  
$$T=-\lim_{q\to +\infty}{\frac{p^*(q)+q^2}{q}}<\infty,$$
solutions of (\ref{x}) exist on $[0, +\infty)$.
\\
Since a solution $x(t)=\Theta_{\tau}(t)$ of (\ref{x}) equals  $-\Theta_{-\tau}(-t)$ if 
$t\le 0$, solutions of (\ref{x}) exist on $(-\infty, +\infty)$.
\hfill  $\Box$ 
\section{\bf Asymptotic  behaviour  at infinity}
\noindent
In this chapter we deal with the properties of 
the solutions of (\ref{u}) for $\tau\in (-T, T)$
at $r=0$ and $r=+\infty$.
\\
In this chapter initial value problem
(\ref{u}) is oft replaced by the autonomous initial value problem
(\ref{x}) and another initial value problem which 
is more convenient when $r\to 0$ or $r\to +\infty$.
\begin{theo}
There is a constant $T$ such that for any solution $\Psi_{\tau}(r)$
of (\ref{u}) with $\tau\in (-T, T)$ it holds
$$\Psi_{\tau}'(r)=\frac{1}{r^2}\xi_{\tau}(r^2)=
\zeta_{\tau} (\frac{1}{r^2}),$$
$$\left(\Psi_{\tau }''r^2 + \Psi_{\tau }'r - \Psi_{\tau }\right)
\Psi_{\tau }'^2r^2=\mu_{\tau}(r^2)r=
\nu_{\tau}(\frac{1}{r^2})\frac{1}{r},$$
where the functions $\xi_{\tau}$, $\zeta_{\tau}$,
$\nu_{\tau}$, $\mu_{\tau}$ are of class 
$C^{\infty}$ and 
$\xi_{\tau}(0)\ne 0,\  \zeta_{\tau}(0)\ne 0$.
\label{thsm}
\end{theo}
\noindent
{\em Proof.} Initial value problem (\ref{u}) 
can be replaced by another  initial value problem.
The change of variables $r\to s=r^{-2}$, $u\to g=r^{-1}u$
changes (\ref{u}) to the following initial value problem
\begin{equation}
g'''=3\frac{g'' g'}{g-2g's} + s\frac{g''^2}{g-2g's},
 \  g(1)=0, g'(1)=-\frac{1}{2}, g''(1)=\tau
\label{g}
\end{equation}
We will prove that solutions of (\ref{g})
for parameter $\tau\in (-T, T)$
satisfy the regular initial conditions at $s=0$:
$g(0)\ne 0, g'(0)<\infty, 0<g''(0)<\infty$.
\medbreak
\noindent
Compute now for a solution of (\ref{x}):
\begin{equation}
x'(t)=(\exp t)(g-2g's),
\label{d}
\end{equation}
\begin{equation}
x''(t)-x(t)=4\exp {(-3t)}g'' (s),
\label{dd}
\end{equation}
remember that $s=\exp{(-2t)}$.
\medbreak
\noindent
Let us consider system (\ref{syst1}).
Since the eigenvalues of  the corresponding linear system are equal to $-2$
and $-4$, there exist 
functions $P$ and $Q$ of class $C^1$ such that 
 $$q=Q(s)=1+Cs+o_1(s),$$ $$p=P(s)=-2Cs+o_2(s),$$
where $C$ is a constant, see \cite{Har}.
\\
We may write:
$$R(t)=\int^{t}q=t+\psi_1(s).$$
Let us show that $\psi_1$ is of class  $C^1$. By differentiation we obtain
$$\frac{d\psi_1(s)}{ds}=-\frac{q-1}{2s}=-\frac{Cs+o_1(s)}{2s}\in C^{0}.$$
Write now for a solution of (\ref{x})
$$x(t)=\exp {R(t)}=(\exp t)\exp {\psi_1(s)}=(\exp t)g(s),$$
where we have used that $g=r^{-1}u=\exp(-t)u=\exp(-t)x(t)$,
and, therefore, $g(0)=\exp {\psi_1(0)}\ne 0$. 
\medbreak
\noindent
For a solution of (\ref{x}) which can be extended into infinity it
 holds  $$x'(t)=x(t)q(t)=(\exp t)(\exp {\psi_1(s)})(1+Cs+o_1(s))= $$ $$
=(\exp t) g(s)(1+Cs+o_1(s)), s=\exp (-2t).$$ 
On the other hand, from (\ref{d}) we obtain 
$$g(s)(1+Cs+o_1(s))=g(s)-2sg'(s). $$ 
Thus, for any solution of  (\ref{x}) where the corresponding orbit of (\ref{syst1}) converges to the 
singular point $q=1, p=0$ it holds $g'(0)<\infty$.
\medbreak
\noindent
Let us show that for these solutions  it holds $0<g''(0)<\infty$.
\\
We compute
$$x''(t)-x(t)=x(t)(p+q^2-1)=(\exp t)(\exp {\psi_1(s)})(P(s)+Q^2(s)-1)=$$
$$=(\exp t)(\exp {\psi_1(s)})(-2Cs+o_2(s)+(1+Cs+o_1(s))^2-1).$$
Thus,
$$x''(t)-x(t)=(\exp t)\psi_2(s),$$ where
$\psi_2\in C^1$, $\psi_2(s)=o_3(s)$.
Rewrite the differential equation in (\ref{x}) in the following form
$$(x''' - x')x'=-3(x'' - x)x - 2(x'' -x)^2.$$
We obtain  either $x''(t)-x(t)\equiv 0$   
or
$$(\log (x''(t)-x(t)))'=-3\frac{1}{q}-2(x''(t)-x(t))(\exp t\exp {\psi_1(s)}(1+Cs+o_1(s)))^{-1}=$$
$$=-3\frac{1}{q}-2{\psi_2(s)}
(\exp {\psi_1(s)}(1+Cs+o_1(s)))^{-1}.$$ 
So, $$(\log (x''(t)-x(t)))'= -3(1-Cs)+o_4(s).$$ 
Then it follows 
$$\left(\log \left((\exp t)\psi_2(s)\right)\right)'=-3(1-Cs)+o_4(s).$$ Thus,
$$(t+\log \psi_2(s))'=-3(1-Cs)+o_4(s).$$
 By integrating we get
$$(\log \psi_2(s))=-4t-\frac{3C}{2}s+\psi_3(s),$$ where
$\psi_3\in C^1$.
\\
Thus,
$$ \psi_2(s)=(\exp(-4t))\exp ({-\frac{3C}{2}s+\psi_3(s)})=
s^2\psi_4(s),$$ since $s=\exp(-2t)$. So,
$\psi_4 \in C^1$ and  $\psi_4(0)=const \ne 0$. Taking into account
(\ref{dd}), we obtain
$g''(0)=\frac{1}{4}\psi_4(0)$ and therefore $0<g''(0)<\infty$.
\\
So, we conclude that 
the solutions of (\ref{g}) for $\tau\in (-T, T)$ 
are of class  $C^{\infty}$ in zero.
\medbreak
\noindent
Let $u(r)=\Psi_{\tau}(r)$ and $x(t)=\Theta_{\tau}(t)$ be solutions of (\ref{u}) and (\ref{x})
respectively.
Then we get that
$\Psi_{\tau}(r)=\Theta_{\tau}(\log r)$ and 
$$\Psi'_{\tau}(r)=\Theta'(\log r)
\frac{1}{r}=g(\frac{1}{r^2})- 2g'(\frac{1}{r^2})\frac{1}{r^2}=
\zeta_{\tau}(\frac{1}{r^2}),$$ 
where $\zeta_{\tau}$ is of class $C^{\infty}$.
From  (\ref{d}) and the condition $g(0)\ne 0$ it follows
that  $\zeta_{\tau}(0)\ne 0$.
\\
Using (\ref{d}) and (\ref{dd}) we may compute
$$\left(\Psi_{\tau }''r^2 + \Psi_{\tau }'r - \Psi_{\tau }\right)
\Psi_{\tau }'^2r^2=\left(\Theta_{\tau}''(t) - 
\Theta_{\tau}(t) \right)\Theta_{\tau}'^2(t)=$$
 $$4\exp(-t)g''(s)(g-2g'(s)s)^2=
4\frac{1}{r}\left(g(\frac{1}{r^2}) - 
2g'(\frac{1}{r^2})\frac{1}{r^2} \right)^2g''(\frac{1}{r})=
\nu_{\tau}(\frac{1}{r^2})\frac{1}{r},$$
where $\nu_{\tau}$ is of class $C^{\infty}$.
\\
Since a solution $x(t)=\Theta_{\tau}(t)$ of (\ref{x}) equals  $-\Theta_{-\tau}(-t)$ if 
$t\le 0$,
in the same way, we obtain
$$\Psi'_{\tau}(r)=\frac{1}{r^2}\xi_{\tau}(r^2),$$
where $\xi_{\tau}\in C^{\infty}$, $\xi_{\tau}(0)\ne 0$ and 
$$\left(\Psi_{\tau }''r^2 + \Psi_{\tau }'r - \Psi_{\tau }\right)
\Psi_{\tau }'^2r^2=\mu_{\tau}(r^2)r,$$
where $\mu_{\tau}$ is of class $C^{\infty}$.
\hfill  $\Box$ 
\section{\bf Smooth conservative systems on $S^2$}
\noindent
In this chapter we propose new examples of complete integrable conservative 
systems on $S^2$ in terms of  solutions of (\ref{u}).
\begin{theo} The one-parameter family of conservative systems with energy (\ref{ms}),
where $\Psi\in \{\Psi_{\tau}\}_{\tau \in (0, T)}$ 
and $\Psi_{\tau}$ is a solution of (\ref{u}), is a one-parameter 
family of smooth systems on $S^2$
possessing  an  integral  cubic  in momenta which is nontrivial, i.e.
there is no quadratic or linear integral.
\\
 These systems cannot be 
obtained one from another by a change of variables.
\\
This family of integrable conservative systems on $S^2$   does not include 
the case of
Goryachev-Chaplygin.
\end{theo}
\noindent
{\em Proof.} 
Prove now that the kinetic energy of (\ref{ms}) for any $\tau\in (0, T)$
 is a $C^{\infty}$ metric
on $S^2$. Using theorem \ref{thsm}, we may write the kinetic energy
as
$$\frac{r^2d\varphi^2 + dr^2}{\xi_{\tau}(r^2)}=
\frac{{\tilde r}^2d{\tilde \varphi}^2 + d{\tilde r}^2}
{\zeta_{\tau}({\tilde r}^2)},$$
where $\tilde r=\frac{1}{r}$, $\tilde \varphi=-\varphi$.
Since functions $\xi_{\tau}$ and $\zeta_{\tau}$ are of class $C^{\infty}$
 and, moreover,
$\xi_{\tau}(0)\ne 0$, $\zeta_{\tau}(0)\ne 0$, it follows that
 the kinetic energy is a smooth metric on $S^2$.
\\
The potential of   (\ref{ms}) also can be written in polar coordinates:
$$V_{\tau}=\mu_{\tau}(r^2)r\cos \varphi=\nu_{\tau}({\tilde r}^2)\tilde 
r\cos \tilde \varphi.$$
Thus, the potential is also a smooth function in 
$r\cos \varphi$, $r\sin \varphi$, $0\le r <\infty$ and 
in ${\tilde r}\cos \tilde \varphi$, ${\tilde r}\sin \tilde \varphi$,
 $0\le {\tilde r}<\infty$ and therefore,
 a smooth function on $S^2$. 
So, any system with energy  (\ref{ms}) is 
a $C^{\infty}$ conservative system  on $S^2$.
\medbreak
\noindent
We prove now the integrability of these systems.
\\
As mentioned above, if  $\Psi_{\tau}$ satisfies (\ref{u})
 then the function
$\Theta_{\tau}=\Psi_{\tau}\circ\exp$ satisfies (\ref{x}).
\\ 
In the coordinate system $y=\log r$, $\varphi=\varphi$ energy 
(\ref{ms}) has  the following form:
\begin{equation}
H_{\tau}={\Theta_{\tau}}'^2(y)(p_{\varphi}^2+p_y^2) - {\Theta_{\tau}}'^2(y)({\Theta_{\tau}}''(y)-{\Theta_{\tau}}(y))\cos {\varphi}, 
\label{x_h}
\end{equation}
where $\Theta_{\tau}\in \{\Theta_{\tau}\}_{\tau \in (0, T)}$.
\\
We will prove  that the following  polynomial cubic in momenta
$$
F_{\tau}=p_{\varphi}^3 + \frac{3}{2}\left
(\Theta_{\tau}(y)\cos {\varphi} {p_{\varphi}}
 - \Theta_{\tau}'(y)\sin {\varphi} {p_y}\right)
=p_{\varphi}^3 + E_{\tau}({\varphi},y,p_{\varphi},p_y)
$$
is an integral 
of the system with the Hamiltonian 
(\ref{x_h}), i.e. $\{F_{\tau}, H_{\tau}\}\equiv 0$.
\\
Denote  by $\hat H_{\tau}= {\Theta_{\tau}}'^2(y)(p_{\varphi}^2+p_y^2)$
 the Hamiltonian of 
the geodesic flow of  metric
$$\hat {ds_{\tau}}^2=\frac{1}{{\Theta_{\tau}}'^2(y)}
\left( d{\varphi}^2 + dy^2\right).$$
Write now  energy as  $ H_{\tau}=\hat H_{\tau}+V_{\tau}$ where 
the potential $V_{\tau}$ has the following form
\begin{equation}
 V_{\tau}= - {\Theta_{\tau}}'^2(y)({\Theta_{\tau}}''(y)-
{\Theta_{\tau}}(y))\cos {\varphi}.
\label{pot}
\end{equation}
Note that $\{\hat H, p_{\varphi}^3\}\equiv 0$. 
In order to prove that $\{F_{\tau}, H_{\tau}\}\equiv 0$ 
we must prove that $\{V_{\tau}, p_{\varphi}^3\} + \{\hat H_{\tau}, 
 E_{\tau}({\varphi},y,p_{\varphi},p_y)\}\equiv 0$ and $\{V_{\tau},
  E_{\tau}({\varphi},y,p_{\varphi},p_y)\}\equiv 0$.
\\
Compute $$\{V_{\tau}, p_{\varphi}^3\} + \{\hat H_{\tau}, 
 E_{\tau}({\varphi},y,p_{\varphi},p_y)\}=-3{\Theta_{\tau}}'^2({\Theta_{\tau}}'' 
-{\Theta_{\tau}})\sin {\varphi} p_{\varphi}^2 - 
3{\Theta_{\tau}}'^2({\Theta_{\tau}}\sin {\varphi} p_{\varphi} + $$ $$
{\Theta_{\tau}}'\cos {\varphi} p_y)p_{\varphi}  + 
3{\Theta_{\tau}}'^2({\Theta_{\tau}}'\cos {\varphi} p_{\varphi} -
 {\Theta_{\tau}}''\sin {\varphi} p_y)p_y + 
3{\Theta_{\tau}}'^2{\Theta_{\tau}}''\sin {\varphi} (p_{\varphi}^2 + p_y^2)\equiv 0.$$
\\
By a computation we obtain
$\{V_{\tau}, E_{\tau}({\varphi},y,p_{\varphi},p_y)\}=$
$$-\frac{3}{2}\cos {\varphi} \sin {\varphi} \Theta_{\tau}'^2\left
(\Theta_{\tau}'\Theta_{\tau}''' -  \Theta_{\tau}^2 - 
\Theta_{\tau}''\Theta_{\tau} - \Theta_{\tau}'^2 + 2\Theta_{\tau}''^2\right)\equiv 0.$$
\\
Let us assume that  a system of this family has an integral
 which is independent
of  the energy and which is a polynomial of second degree with
respect to momenta (clearly, this assumption includes
 the case of linear integrals). So, there is an  integral $\tilde F_{\tau}$ of (\ref{x_h}) which 
is quadratic in momenta.
Thus, $\tilde F_{\tau}= A_{\tau}(p_{\varphi}, p_{y}, \varphi, y) + B_{\tau}(\varphi, y)$ 
where $ A_{\tau}(p_{\varphi}, p_{y}, \varphi, y)$
is a polynomial of  second degree with respect to momenta.
We may write $\{\tilde F_{\tau}, H_{\tau}\}=\{
A_{\tau}(p_{\varphi}, p_{y}, \varphi, y) + B_{\tau}(\varphi, y), 
\hat H_{\tau} +V_{\tau}\}\equiv 0$
and, therefore, $\{A_{\tau}(p_{\varphi}, p_{y}, \varphi, y), \hat H_{\tau}\}\equiv 0$. 
\\
Thus, the geodesic flow of $\hat {ds_{\tau}}^2$ has an integral which is
a polynomial of  second degree with respect to momenta
and, therefore, in some Liouville coordinate system $u, v$ it holds
\begin{equation}
\hat H_{\tau}=\frac{1}{\alpha_1(u)+\alpha_2(v)}(p_u^2+p_v^2), V_{\tau}=\frac{\alpha_3(u)
+\alpha_4(v)}{\alpha_1(u)+\alpha_2(v)}
\label{lio}
\end{equation}
where $\alpha_1$, $\alpha_2$, $\alpha_3$, $\alpha_4$ are some functions.
\\
Geodesic flows on $S^2$ 
 possesing an  integral  quadratic in momenta have been described in \cite{kol}.
 In particular, it has been shown
that Liouville coordinate systems for all metrics on  $S^2$   of
nonconstant curvature are unique up to linear combinations.
Noting that (\ref{pot}) for any $\tau\in (0, T)$ has not
form (\ref{lio}) in any coordinate systems 
 $ay+b$, $a{\varphi} +d$, where $a$, $b$, $d$ are some constants,
we prove that no  system from our series possesses an 
integral quadratic in momenta. 
Moreover, since the solutions of (\ref{x}) for $\tau\in (0, T)$
 cannot  be obtained one from another by a linear change of variables,
systems with the energy (\ref{ms}) also cannot be obtained
 one from another by a change of variables. 
\medbreak
\noindent
Prove  that our family of smooth conservative systems on $S^2$ does not 
include the  well-known  case of Goryachev-Chaplygin from the dynamics 
of a rigid body.
\\
In the case of Goryachev-Chaplygin there is also a coordinate system $\varphi, y$
where  energy has the following form
$$H=\gamma(y)(p_{\varphi}^2 + p_y^2) + \beta(y)\cos \varphi$$
for some functions $\gamma$ and $\beta$, see  \cite{BKF}.
Since $\hat H=\gamma(y)(p_{\varphi}^2 + p_y^2) $ is 
the Hamiltonian of the geodesic flow of 
a smooth metric on $S^2$,
 this coordinate system is unique up to linear combinations.
\\
On the other hand there is an integral $F=p_{\varphi}^3 + \kappa p_{\varphi}\hat H + 
G(p_{\varphi}, p_{y}, \varphi, y)$ where constant 
$\kappa \ne 0$ and $G(p_{\varphi}, p_{y}, \varphi, y)$ is a polynomial linear in momenta
(A.V. Bolsinov, private communication), see also \cite{BKF}.
\\
If our series includes the case of Goryachev-Chaplygin then it follows that in the case
 of Gorychev-Chaplygin
there are two independent first integrals. 
So, our series of integrable conservative systems on $S^2$ does not 
include  the case of Goryachev-Chaplygin.
\hfill  $\Box$ 
\medbreak
\noindent
{\bf Commentar.}
Numerically we got that  the value of $T$ is equal to $\approx $0,57735.
\medbreak
\noindent
{\bf Acknowledgement.} The author wishes to thank  
 A.V. Bolsinov,  A.T. Fomenko and A.M. Stepin for their
interest in her work.
\\
The author would like to thank Prof. G. Huisken and Arbeitsbereich 
"Analysis" of Eberhard-Karls-Universit\"at T\"ubingen 
 for their hospitality. 

\bigbreak
\small  {\em Mathematisches Institut, Universit\"at T\"ubingen} \\   
 \small  {\em Auf der Morgenstelle 10,  72076  T\"ubingen, Germany} \\
 \small  {\em e-mail: \ \ lena@moebius.mathematik.uni-tuebingen.de} 
\bigbreak
\small  {\em Chair of Geometry, Nizhny Novgorod State 
Pedagogical University} \\
\small  {\em 603000 Russia, Nizhny Novgorod, ul. Ulyanova 1}
\end{document}